\documentclass[a4paper,12pt,final]{amsart}
\usepackage{times,a4wide,mathrsfs,amssymb,amsmath,amsthm,enumerate,xypic,tikzsymbols,dsfont}

\newcommand{\C}{\mathbb{C}}

\newcommand{\QQ}{\mathbb{Q}}
\newcommand{\NN}{\mathbb{N}}
\newcommand{\PP}{\mathbb{P}}

\newcommand{\OO}{\mathcal O}

\newcommand{\YY}{\mathcal Y}

\newcommand{\Cc}{\mathcal C}

\newcommand{\MM}{\mathcal M}

\newcommand{\wt}{\widetilde}
\newcommand{\rom}{\romannumeral}

\newcommand{\one}{\mathds{1}}

\DeclareMathOperator{\ima}{Im}

\DeclareMathOperator{\sym}{Sym}
\DeclareMathOperator{\Gr}{Gr}

\DeclareMathOperator{\Dec}{Dec}
\DeclareMathOperator{\Pf}{Pf}

\newtheorem{theorem}{Theorem}[section]

\newtheorem{lemma}[theorem]{Lemma}

\newtheorem{corollary}[theorem]{Corollary}
\newtheorem{proposition}[theorem]{Proposition}

\newtheorem{convention}{Conventions}

\newtheorem{nonumbering}{Theorem}

\newtheorem{nonumberingc}{Corollary}

\theoremstyle{definition}
\newtheorem{remark}[theorem]{Remark}
\newtheorem{definition}[theorem]{Definition}

\newtheorem{notation}[theorem]{Notation}

\newtheorem{nonumberingt}{Acknowledgments}

\begin{document}

\author[Robert Laterveer]
{Robert Laterveer}

\address{Institut de Recherche Math\'ematique Avanc\'ee,
CNRS -- Universit\'e 
de Strasbourg,\
7 Rue Ren\'e Des\-car\-tes, 67084 Strasbourg CEDEX,
FRANCE.}
\email{robert.laterveer@math.unistra.fr}

\title[Some new Fano varieties with MCK]{Some new Fano varieties with a multiplicative Chow--K\"unneth decomposition}

\begin{abstract} Let $Y$ be a smooth dimensionally transverse intersection of the Grassmannian $\Gr(2,n)$ with 3 Pl\"ucker hyperplanes. We show that $Y$
 admits a multiplicative Chow--K\"unneth decomposition, in the sense of Shen--Vial. 
As a consequence, a certain tautological subring of the Chow ring of powers of $Y$ injects into cohomology.
 \end{abstract}

\thanks{Supported by ANR grant ANR-20-CE40-0023.}

%

\maketitle

\section{Introduction}

Given a smooth projective variety $Y$ over $\C$, let $A^i(Y):=CH^i(Y)_{\QQ}$ denote the Chow groups of $Y$, i.e. the groups of codimension $i$ algebraic cycles on $Y$ with $\QQ$-coefficients, modulo rational equivalence. Let us write $A^\ast_{hom}(Y)$ and $A^\ast_{AJ}(Y)$ for the subgroups of homologically trivial (resp. Abel--Jacobi trivial) cycles.
Intersection product defines a ring structure on $A^\ast(Y)=\bigoplus_i A^i(Y)$, the {\em Chow ring\/} of $Y$ \cite{F}. In the case of K3 surfaces, this ring structure has a peculiar property:

\begin{theorem}[Beauville--Voisin \cite{BV}]\label{K3} Let $S$ be a K3 surface. 
The $\QQ$-subalgebra
  \[  R^\ast(S):=  \bigl\langle  A^1(S), c_j(S) \bigr\rangle\ \ \ \subset\ A^\ast(S) \]
  injects into cohomology under the cycle class map.
  \end{theorem}

Inspired by the remarkable behaviour of K3 surfaces and of abelian varieties, Beauville \cite{Beau3} has famously conjectured that for certain special varieties, the Chow ring should admit a {\em multiplicative splitting\/}. To make concrete sense of Beauville's elusive ``splitting property conjecture'', Shen--Vial \cite{SV} have introduced the concept of {\em multiplicative Chow--K\"unneth decomposition\/}; 
let us abbreviate this to ``MCK decomposition''.

What can one say about the class of special varieties admitting an MCK decomposition ? This class is not yet well-understood. Varieties with $A^\ast_{hom}()=0$ (i.e. varieties with {\em trivial Chow groups\/}) admit an MCK decomposition, for trivial reasons. The question becomes interesting for varieties with $A^\ast_{AJ}()=0$ (conjecturally, these are exactly the varieties with Hodge level at most 1, i.e. the Hodge numbers $h^{p,q}$ are zero for $\vert p-q\vert>1$).
It is known that hyperelliptic curves have an MCK decomposition \cite[Example 8.16]{SV}, but the very general curve of genus $\ge 3$ does not have an MCK decomposition \cite[Example 2.3]{FLV2} (for more details, cf. subsection \ref{ss:mck} below). Also, there exist Fano threefolds that do not admit an MCK decomposition.
On the positive side, here are some higher-dimensional varieties with Hodge level 1 that are known to have an MCK decomposition:

\begin{itemize}

\item cubic threefolds and cubic fivefolds \cite{Diaz}, \cite{FLV2};

\item Fano threefolds of genus 8 \cite{g8};

\item complete intersections of 2 quadrics \cite{2q};

\item Gushel--Mukai fivefolds \cite{46}.
\end{itemize}

The goal of the present note is to add some new varieties with Hodge level 1 to this list:

\begin{nonumbering}[=Theorem \ref{main}] Let $Y$ be a smooth dimensionally transverse intersection
  \[ Y:= \Gr(2,n)\cap H_1\cap H_2 \cap H_3\ \ \ \subset\  \PP^{{n\choose 2}-1}  \ ,\]
  where $\Gr(2,n)$ denotes the Grassmannian of 2-dimensional linear subspaces of a fixed $n$-dimensional vector space, and the $H_j$ are Pl\"ucker hyperplanes.
Then $Y$ has an MCK decomposition. 
\end{nonumbering}

In case $n$ is odd, a variety $Y$ as in Theorem \ref{main} has trivial Chow groups and so the statement is vacuously true. In case $n$ is even, there is a curve $C$ naturally associated to $Y$, and one has a relation of Chow motives
  \begin{equation}\label{isomot}  h(Y) \cong h(C)((1-\dim Y)/2)\oplus \bigoplus \one(\ast)\ \ \ \hbox{in}\ \MM_{\rm rat} \end{equation}
  (cf. Theorem \ref{main0}).
The relation between $Y$ and $C$ has previously been studied on the level of Hodge theory in \cite{Don}, and on the level of derived categories in \cite{Kuz-1}, \cite{Kuz0}. As a result of independent interest, we prove here (Theorem \ref{main0}) that the relation \eqref{isomot} also holds on the level of Chow motives.

The existence of an MCK decomposition has profound intersection-theoretic consequences. This is exemplified by the following corollary, which is about a certain {\em tautological subring\/} of the Chow ring of powers of $Y$:

\begin{nonumberingc}[=Corollary \ref{cor1}]
Let $Y$ be as in Theorem \ref{main}, and $m\in\NN$. Let
  \[ R^\ast(Y^m):=\Bigl\langle (p_i)^\ast \ima\bigl(  A^\ast(\Gr(2,n))\to A^\ast(Y)\bigr)   \, ,\  (p_{ij})^\ast(\Delta_Y)\Bigr\rangle\ \subset\ \ \ A^\ast(Y^m)   \]
  be the $\QQ$-subalgebra generated by (pullbacks of) cycles coming from the Grassmannian and the diagonal $\Delta_Y\in A^\ast(Y\times Y)$. (Here $p_i$ and $p_{ij}$ denote the various projections from $Y^m$ to $Y$ resp. to $Y\times Y$).
  The cycle class map induces injections
   \[ R^\ast(Y^m)\ \hookrightarrow\ H^\ast(Y^m,\QQ)\ \ \ \hbox{for\ all\ }m\in\NN\ .\]
   \end{nonumberingc}

Corollary \ref{cor1} is somewhat surprising, because the corresponding statement for the associated curve $C$ is {\em false\/}: in general it is {\em not\/} true that the $\QQ$-subalgebra
  \[ \Bigl\langle  (p_i)^\ast (h),   \, (p_{ij})^\ast(\Delta_C)\Bigr\rangle\ \subset\ \ \ A^\ast(C^m)   \]
  injects into cohomology (cf. Proposition \ref{neg} for the precise statement). This means that the injection
  \[ A^\ast(C^m)\ \hookrightarrow\ A^\ast(Y^m)  \]
  induced by \eqref{isomot} does {\em not\/} send tautological cycles to tautological cycles !

  Let us end this introduction with an open question. In view of Theorem \ref{main}, one might ask whether more generally smooth complete intersections of Grassmannians $\Gr(k,n)$ with an {
  arbitrary} number of Pl\"ucker hyperplanes have an MCK decomposition. This concerns in particular the Debarre--Voisin 20folds
  \[ \Gr(3,10)\cap H\ , \]
  which are Fano varieties of K3 type \cite{DV}, and also the
Fano eightfolds
    \[ \Gr(2,8)\cap H_1\cap\cdots\cap H_4\ ,\] 
 which are again of K3 type \cite{ST}, \cite{FM}. Such varieties (being of Hodge level $>1$) are out of scope of the argument of the present note.
   
 \vskip0.5cm

\begin{convention} In this note, the word {\sl variety\/} will refer to a reduced irreducible scheme of finite type over $\C$. A {\sl subvariety\/} is a (possibly reducible) reduced subscheme which is equidimensional. 

{\bf All Chow groups will be with rational coefficients}: we denote by $A_j(Y)$ the Chow group of $j$-dimensional cycles on $Y$ with $\QQ$-coefficients; for $Y$ smooth of dimension $n$ the notations $A_j(Y)$ and $A^{n-j}(Y)$ are used interchangeably. 
The notations $A^j_{\rm hom}(Y)$ and $A^j_{\rm AJ}(Y)$ will be used to indicate the subgroup of homologically trivial (resp. Abel--Jacobi trivial) cycles.

The contravariant category of Chow motives (i.e., pure motives with respect to rational equivalence as in \cite{Sc}, \cite{MNP}) will be denoted 
$\MM_{\rm rat}$.
\end{convention}

\section{Preliminaries}

\subsection{MCK decomposition}
\label{ss:mck}

\begin{definition}[Murre \cite{Mur}] Let $X$ be a smooth projective variety of dimension $n$. We say that $X$ has a {\em CK decomposition\/} if there exists a decomposition of the diagonal
   \[ \Delta_X= \pi^0_X+ \pi^1_X+\cdots +\pi_X^{2n}\ \ \ \hbox{in}\ A^n(X\times X)\ ,\]
  such that the $\pi^i_X$ are mutually orthogonal idempotents and $(\pi_X^i)_\ast H^\ast(X,\QQ)= H^i(X,\QQ)$.
  
  (NB: ``CK decomposition'' is shorthand for ``Chow--K\"unneth decomposition''.)
\end{definition}

\begin{remark} The existence of a CK decomposition for any smooth projective variety is part of Murre's conjectures \cite{Mur}, \cite{J4}. 
\end{remark}

\begin{definition}[Shen--Vial \cite{SV}] Let $X$ be a smooth projective variety of dimension $n$. Let $\Delta_X^{\rm sm}\in A^{2n}(X\times X\times X)$ be the class of the small diagonal
  \[ \Delta_X^{\rm sm}:=\bigl\{ (x,x,x)\ \vert\ x\in X\bigr\}\ \subset\ X\times X\times X\ .\]
  An {\em MCK decomposition\/} is a CK decomposition $\{\pi_X^i\}$ of $X$ that is {\em multiplicative\/}, i.e. it satisfies
  \[ \pi_X^k\circ \Delta_X^{\rm sm}\circ (\pi_X^i\times \pi_X^j)=0\ \ \ \hbox{in}\ A^{2n}(X\times X\times X)\ \ \ \hbox{for\ all\ }i+j\not=k\ .\]
  
 (NB: ``MCK decomposition'' is shorthand for ``multiplicative Chow--K\"unneth decomposition''.) 
  \end{definition}
  
  \begin{remark} The small diagonal (seen as a correspondence from $X\times X$ to $X$) induces the {\em multiplication morphism\/}
    \[ \Delta_X^{\rm sm}\colon\ \  h(X)\otimes h(X)\ \to\ h(X)\ \ \ \hbox{in}\ \MM_{\rm rat}\ .\]
 Let us assume $X$ has a CK decomposition
  \[ h(X)=\bigoplus_{i=0}^{2n} h^i(X)\ \ \ \hbox{in}\ \MM_{\rm rat}\ .\]
  By definition, this decomposition is multiplicative if for any $i,j$ the composition
  \[ h^i(X)\otimes h^j(X)\ \to\ h(X)\otimes h(X)\ \xrightarrow{\Delta_X^{\rm sm}}\ h(X)\ \ \ \hbox{in}\ \MM_{\rm rat}\]
  factors through $h^{i+j}(X)$.
  
  If $X$ has an MCK decomposition, then setting
    \[ A^i_{(j)}(X):= (\pi_X^{2i-j})_\ast A^i(X) \ ,\]
    one obtains a bigraded ring structure on the Chow ring: that is, the intersection product sends $A^i_{(j)}(X)\otimes A^{i^\prime}_{(j^\prime)}(X) $ to  $A^{i+i^\prime}_{(j+j^\prime)}(X)$.
    
   It is expected that for any $X$ with an MCK decomposition, one has
    \[ A^i_{(j)}(X)\stackrel{??}{=}0\ \ \ \hbox{for}\ j<0\ ,\ \ \ A^i_{(0)}(X)\cap A^i_{\rm hom}(X)\stackrel{??}{=}0\ ;\]
    this is related to Murre's conjectures B and D, that have been formulated for any CK decomposition \cite{Mur}.

  The property of having an MCK decomposition is restrictive, and is closely related to Beauville's ``splitting property conjecture'' \cite{Beau3}. 
  To give an idea: hyperelliptic curves have an MCK decomposition \cite[Example 8.16]{SV}, but the very general curve of genus $\ge 3$ does not have an MCK decomposition \cite[Example 2.3]{FLV2}. As for surfaces: a smooth quartic in $\PP^3$ has an MCK decomposition, but a very general surface of degree $ \ge 7$ in $\PP^3$ should not have an MCK decomposition \cite[Proposition 3.4]{FLV2}. There are examples of Fano threefolds that do not admit an MCK decomposition \cite[Example 1.11]{FLV2}.
  
For more detailed discussion, and examples of varieties with an MCK decomposition, we refer to \cite[Section 8]{SV}, as well as \cite{V6}, \cite{SV2}, \cite{FTV}, \cite{37}, \cite{38}, \cite{39}, \cite{40}, \cite{2q}, \cite{3q}, \cite{FLV2}.
   \end{remark}

 \subsection{The Franchetta property}
 \label{ss:fr}

 \begin{definition} Let $\YY\to B$ be a smooth projective morphism, where $\YY, B$ are smooth quasi-projective varieties. We say that $\YY\to B$ has the {\em Franchetta property in codimension $j$\/} if the following holds: for every $\Gamma\in A^j(\YY)$ such that the restriction $\Gamma\vert_{Y_b}$ is homologically trivial for the very general $b\in B$, the restriction $\Gamma\vert_b$ is zero in $A^j(Y_b)$ for all $b\in B$.
 
 We say that $\YY\to B$ has the {\em Franchetta property\/} if $\YY\to B$ has the Franchetta property in codimension $j$ for all $j$.
 \end{definition}
 
 This property is studied in \cite{PSY}, \cite{BL}, \cite{FLV}, \cite{FLV3}.
 
 \begin{definition} Given a family $\YY\to B$ as above, with $Y:=Y_b$ a fiber, we write
   \[ GDA^j_B(Y):=\ima\Bigl( A^j(\YY)\to A^j(Y)\Bigr) \]
   for the subgroup of {\em generically defined cycles}. 
  In a context where it is clear to which family we are referring, the index $B$ will often be suppressed from the notation.
  \end{definition}
  
  With this notation, the Franchetta property amounts to saying that $GDA^\ast_B(Y)$ injects into cohomology, under the cycle class map, for every fiber $Y$.
  
  There is some flexibility with respect to the base $B$:
  
  \begin{lemma}\label{spread} Let $\YY\to B$ be a smooth projective family, and $B_0\subset B$ the intersection of a countable number of dense open subsets. Then $\YY\to B$ has the Franchetta property if and only if $\YY\to B_0$ has the Franchetta property.
  \end{lemma}
  
  \begin{proof} This follows from a well-known spread lemma \cite[Lemma 3.2]{Vo}.
    \end{proof}

   \subsection{A Franchetta-type result}
   
  \begin{proposition}\label{Frtype} Let $M$ be a smooth projective variety with trivial Chow groups. Let $L_1,\ldots,L_r\to {M}$ be very ample line bundles, and let
  $\YY\to B$ be the universal family of smooth dimensionally transverse complete intersections of type 
    \[ Y={M}\cap H_1\cap\cdots\cap H_r\ ,\ \ \  H_j\in\vert L_j\vert\ .\]
  Assume the fibers $Y=Y_b$ have $H^{\dim Y}_{\rm tr}(Y,\QQ)\not=0$.
  There is an inclusion
    \[ \ker \Bigl( GDA^{\dim Y}_B(Y\times Y)\to H^{2\dim Y}(Y\times Y,\QQ)\Bigr)\ \ \subset\ \Bigl\langle (p_1)^\ast GDA^\ast_B(Y), (p_2)^\ast GDA^\ast_B(Y)  \Bigr\rangle\ .\]
   \end{proposition}
   
   \begin{proof} This is essentially Voisin's ``spread'' result \cite[Proposition 1.6]{V1} (cf. also \cite[Proposition 5.1]{LNP} for a reformulation of Voisin's result). We give a proof which is somewhat different from \cite{V1}. Let $\bar{B}:=\PP H^0({M},L_1\oplus\cdots\oplus L_r)$ (so $B\subset \bar{B}$ is a Zariski open), and let us consider the projection
   \[ \pi\colon\ \  \YY\times_{\bar{B}} \YY\ \to\ M\times M\ .\]
   Using the very ampleness assumption, one finds that $\pi$ is a $\PP^s$-bundle over $({M}\times {M})\setminus \Delta_{{M}}$, and a $\PP^t$-bundle over $\Delta_{{M}}$.
   That is, $\pi$ is what is termed a {\em stratified projective bundle\/} in \cite{FLV}. As such, \cite[Proposition 5.2]{FLV} implies the equality
      \begin{equation}\label{stra} GDA^\ast_B(Y\times Y)= \ima\Bigl( A^\ast(M\times M)\to A^\ast(Y\times Y)\Bigr) +  \Delta_\ast GDA^\ast_B(Y)\ ,\end{equation}
      where $\Delta\colon Y\to Y\times Y$ is the inclusion along the diagonal. As $M$ has trivial Chow groups, $A^\ast(M\times M)$ is generated by $A^\ast(M)\otimes A^\ast(M)$. 
      Base-point freeness of the $L_j$ implies that 
        \[  GDA^\ast_B(Y)=\ima\bigl( A^\ast(M)\to A^\ast(Y)\bigr)\ .\]
       The equality \eqref{stra} thus reduces to
      \[ GDA^\ast_B(Y\times Y)=\Bigl\langle (p_1)^\ast GDA^\ast_B(Y), (p_2)^\ast GDA^\ast_B(Y), \Delta_Y\Bigr\rangle\ \]   
      (where $p_1, p_2$ denote the projection from $S\times S$ to first resp. second factor). The assumption that $Y$ has non-zero transcendental cohomology
      implies that the class of $\Delta_Y$ is not decomposable in cohomology. It follows that
      \[ \begin{split}  \ima \Bigl( GDA^{\dim Y}_B(Y\times Y)\to H^{2\dim Y}(Y\times Y,\QQ)\Bigr) =&\\
       \ima\Bigl(  \Dec^{\dim Y}(Y\times Y)\to H^{2\dim Y}(Y\times Y,\QQ)\Bigr)& \oplus \QQ[\Delta_Y]\ ,\\
       \end{split}\]
      where we use the shorthand 
       \[ \Dec^j(Y\times Y):= \Bigl\langle (p_1)^\ast GDA^\ast_B(Y), (p_2)^\ast GDA^\ast_B(Y)\Bigr\rangle\cap A^j(Y\times Y) \ \]     
       for the {\em decomposable cycles\/}. 
       We now see that if $\Gamma\in GDA^{\dim Y}(Y\times Y)$ is homologically trivial, then $\Gamma$ does not involve the diagonal and so $\Gamma\in \Dec^{\dim Y}(Y\times Y)$.
       This proves the proposition.
         \end{proof}
  
   
\begin{corollary}\label{Frtype2} Let $\YY\to B$ be as in Proposition \ref{Frtype}. Assume that $\YY\to B$ has the Franchetta property. Then for any fiber $Y$ the cycle class map induces an injection
  \[ GDA^{\dim Y}(Y\times Y)\ \hookrightarrow\ H^{2\dim Y}(Y\times Y,\QQ) \ .\]
\end{corollary}

\begin{proof} This is immediate from Proposition \ref{Frtype}: the Franchetta property for $\YY\to B$, combined with the K\"unneth decomposition in cohomology, implies that the right-hand side of Proposition \ref{Frtype} injects into cohomology.
\end{proof}
   
 \subsection{A CK decomposition}

\begin{lemma}\label{ck} Let $M$ be a smooth projective variety with trivial Chow groups. Let $Y\subset M$ be a smooth complete intersection of dimension $\dim Y=d$ defined by ample line bundles.
 The variety $Y$ has a self-dual CK decomposition $\{\pi^j_Y\}$ with the property that
   \[  h^j(Y):=(Y,\pi^j_Y,0) =\oplus \one(\ast)\ \ \ \hbox{in}\ \MM_{\rm rat}\ \ \ \forall\ j\not=d \ .\]
   
   Moreover, this CK decomposition is {\em generically defined\/}: writing $\YY\to B$ for the universal family (of complete intersections of the type of $Y$),
 there exist relative projectors $\pi^j_\YY\in A^{d}(\YY\times_B \YY)$ such that $\pi^j_Y=\pi^j_\YY\vert_{b}$ (where $Y=Y_b$ for $b\in B$). 
     \end{lemma} 
 
 \begin{proof} This is a standard construction, one can look for instance at \cite{Pet} (in case $d$ is odd, which will be the case in this note,  the ``variable motive'' $h(Y)^{\rm var}$ of 
 \cite[Theorem 4.4]{Pet}  coincides with $h^{d}(Y)$).
\end{proof}

\section{Main results}

\subsection{An isomorphism of motives}

\begin{definition}\label{dual} Let $V$ be a vector space of dimension $n$, and let
  \[ \Gr(2,n):= \Gr(2,V)\ \ \subset\ \PP(\wedge^2 V) \]
  be the Grassmannian (parametrizing $2$-dimensional subspaces of $V$) in its Pl\"ucker embedding. Assuming $n$ is even, let
  \[ \Pf\ \subset \ \PP(\wedge^2 V^\vee) \]
  denote the projective dual of $\Gr(2,n)\subset\PP(\wedge^2 V)$, called the {\em Pfaffian.} (The Pfaffian $\Pf$ is a hypersurface of degree $n/2$ and singular locus of codimension $7$.)
  
   Assume $n$ is even. Given a linear subspace $U\subset \wedge^2 V$ of codimension $3$, one can define varieties by intersecting on the Grassmannian side and on the Pfaffian side:
    \[  \begin{split}   Y&=Y_U:=  \Gr(2,V)\cap   \PP(U) \ \ \ \subset\ \PP(\wedge^2 V) \ , \\
                              C&=C_U:= \Pf\,\cap\, \PP(U^\perp)\ \ \ \subset\ \PP(\wedge^2 V^\vee) \ .
                              \end{split}  \]  
                              We say that $Y$ and $C$ are {\em dual\/}.
  For $U$ generic, the intersections $Y$ and $C$ are smooth and dimensionally transverse, of dimension $2(n-2)-3$ resp. 1. 
\end{definition}

\begin{theorem}\label{main0} Let $Y$ be a smooth dimensionally transverse intersection
   \[ Y:=\Gr(2,n)\cap H_1\cap H_2 \cap H_3 \ ,\]
   where the $H_j$ are Pl\"ucker hyperplanes. 
      
   \noindent
   (\rom1) Assume $n$ is odd. Then $A^\ast_{hom}(Y)=0$.
   
   \noindent
   (\rom2) Assume that $n$ is even, and that $Y$ has a smooth dual curve $C$.   
   There is an isomorphism
   \[ h^d(Y)\cong h^1(C)((1-d)/2)\ \ \ \hbox{in}\ \MM_{\rm rat}\ ,\]
   where $d:=\dim Y$ and $h^d(Y)$ is as in Lemma \eqref{ck}.
   \end{theorem}
   
   \begin{proof} This is a special case of \cite[Theorem 3.17]{GrPf}. Since this is crucial to the present note, let us include a (sketch of) proof.
   
    With notation as in Definition \ref{dual}, let us consider 
      \[ Q:= \Bigl\{ (T,\C\omega)\ \in\ \Gr(2,V)\times \PP(U^\perp)\  \big\vert\   \omega\vert_T=0\Bigr\}\ \ \subset\  \Gr(2,V)\times \PP(U^\perp)\  , \]
                  the so-called {\em Cayley hypersurface\/}.  
     There is a diagram
   \begin{equation}\label{diag} \begin{array}[c]{ccccccccc}   && Q_Y & \hookrightarrow & Q & \hookleftarrow & Q_C && \\
                &&&&&&&&\\
                   &{}^{\scriptstyle } \swarrow \ \ && {}^{\scriptstyle p} \swarrow \ \ \ & & \ \ \ \searrow {}^{\scriptstyle q} & & \ \ \searrow {}^{\scriptstyle } & \\
                   &&&&&&&&\\
                   Y & \hookrightarrow & \Gr(2,V) &  &  & & \PP(U^\perp) & \hookleftarrow & C\\
                   \end{array}\end{equation}
                     Here, $C$ is defined to be the empty set for $n$ is odd, and the dual curve $C\subset\Pf$ in case $n$ is even.
                     The morphisms $p$ and $q$ are induced by the natural projections, and the closed subvarieties $Q_Y, Q_C\subset Q$ are defined as $p^{-1}(Y)$ resp. $q^{-1}(C)$. 

The restriction of $p$ to $Q\setminus Q_Y$ is trivial with fibre $Q_u\cong\PP^{1}$, while the restriction of $p$ to $Q_Y$ is Zariski locally trivial with fibre $Q_{Y,y}\cong\PP^{2}$. This allows us to relate the motives of $Q$ and $Y$: an application of the ``motivic Cayley trick'' \cite[Corollary 3.2]{Ji} gives an isomorphism
  \begin{equation}\label{cay} \begin{split}  h(Q)&\ \cong\ h(Y)(-2)\oplus h(\Gr(2,n))\oplus h(\Gr(2,n))(-1)\\
                                                &\ \cong\ h(Y)(-2)\oplus \bigoplus\one(\ast)  \ \ \ \hbox{in}\ \MM_{\rm rat}\ .\\
                                                \end{split}
                                                \end{equation}
  
The restriction of $q$ to $Q_C$ is {\em piecewise trivial\/} (in the sense of \cite[Section 4.2]{Seb}) with constant fiber $F_1$, while the restriction of $q$ to $Q\setminus Q_C$ is piecewise trivial with constant fiber $F_2$. The fibers $F_1$ and $F_2$ are explicitly known; they have only algebraic cohomology \cite[Lemma 3.5]{GrPf}.
This allows to relate $Q$ and $C$ on the level of the Grothendieck ring of varieties, and hence also on the level of cohomology: 
  \begin{equation}\label{right} h(Q) \cong h(C)(2-n)\oplus \bigoplus\one(\ast)\ \ \ \hbox{in}\ \MM_{\rm hom}\ .\end{equation}
  (Here the convention is that $h(C)=0$ in case $n$ is odd.)
  
  Combining \eqref{cay} and \eqref{right}, we find a split injection of homological motives
   \begin{equation}\label{spl}  h^d(Y)\ \hookrightarrow\ h^1(C)((1-d)/2)\ \ \ \hbox{in}\ \MM_{\rm hom}\ .\end{equation}
  
  Let us now consider things family-wise. Writing $B_0\subset \bar{B}:=\PP H^0(\PP(\wedge^2 V),\OO(1)^{\oplus 3})$ for the dense open parametrizing sections such that both $Y_b:=\Gr(2,n)\cap H_1^b\cap H_2^b\cap H_3^b$ and the dual curve $C_b\subset\Pf$ are smooth and dimensionally transverse (and in addition $C_b$ is contained in the non-singular locus $\Pf^\circ\subset\Pf$), we have universal families
    \[  \YY\to B_0\ ,\ \ \Cc\to B_0\ .\] 
 The above construction can be performed for every fiber $Y=Y_b$ of the family $\YY\to B_0$. A Hilbert schemes argument \cite[Proposition 2.11]{GrPf} then allows to find {\em generically defined\/} correspondences (with respect to $B_0$) inducing the split injection \eqref{spl}. Then, the Franchetta-type result (Proposition \ref{Frtype})
 allows to lift the split injection \eqref{spl} to an injection of Chow groups:
      \begin{equation}\label{spl2}  A^\ast_{hom}(Y)= A^\ast\bigl( h^d(Y)\bigr)\ \hookrightarrow\ A^\ast_{hom}\bigl(h^1(C)((1-d)/2)\bigr)=A^1_{hom}(C)\ .\end{equation}
   We conclude from \eqref{spl2} that $A^\ast_{AJ}(Y)=0$ and so
   $Y$ is Kimura finite-dimensional (i.e. $h(Y)$ is finite-dimensional in the sense of \cite{Kim}). Combining \eqref{cay} and \eqref{right}, we find a numerical equality $\dim H^d(Y,\QQ)=\dim H^1(C ,\QQ)$ and so the injection \eqref{spl} is actually an isomorphism of homological motives. Using Kimura finite-dimensionality of both sides, it follows that \eqref{spl} is also
   an isomorphism of Chow motives:
   \[   h^d(Y)\ \xrightarrow{\cong}\ h^1(C)((1-d)/2)\ \ \ \hbox{in}\ \MM_{\rm rat}\ . \]  
   This proves the theorem.
  \end{proof}

 \subsection{Some instances of the Franchetta property}  
 
 \begin{notation}\label{not} Let $\bar{B}$ and $B_0$ be as in the proof of Theorem \ref{main0}, and let $B\supset B_0$ be the set parametrizing smooth dimensionally transverse intersections $Y_b=\Gr(2,n)\cap H_1\cap H_2\cap H_3$; there is a universal family
   \[ \YY\ \to\ B\ .\]
   Assuming $n$ is even, let us write
   \[ \Cc\ \to\ B_0 \]
   for the universal family of smooth dual curves $C_b\subset\Pf^\circ$, as in the proof of Theorem \ref{main0}.
  \end{notation}
  
  \begin{proposition}\label{Fr} The following families have the Franchetta property:
 
 \noindent
 (\rom1) the family $\YY\to B$;
 
 \noindent
 (\rom2) the family $\Cc\to B_0$.
 \end{proposition}
 
 \begin{proof} For (\rom1), let us note that the statement is vacuously true in case $n$ is odd, because then each fiber $Y_b$ has trivial Chow groups (Theorem \ref{main0}(\rom1)). Let us now assume that $n$ is even, say $n=2m$.
We observe that the projection
   \[ \bar{\YY}\ \to\ \Gr(2,n) \]
   is a projective bundle, and so (reasoning with the projective bundle formula, or directly applying \cite[Proposition 5.2]{FLV}) one finds that for any fiber $Y:=Y_b$ there is equality
   \[ GDA^j(Y) = \ima\Bigl(  A^j(\Gr(2,n))\to A^j(Y)\Bigr) \ .\]
   We know from Theorem \ref{main0}(\rom2) that the only non-trivial Chow group of $Y$ is 
   \[  A^{(d+1)/2}(Y)=  A^{(2(n-2)-3+1)/2}(Y)= A^{2m-3}(Y)\ ,\]
   and so we only need to prove that $GDA^{2m-3}(Y)$ injects into cohomology. The Chow ring of the Grassmannian is
   \[ A^\ast(\Gr(2,n))=\bigl\langle h,\, c\bigr\rangle \ ,\]
   where $c:=c_2(Q)\in A^2(\Gr(2,n))$ is the second Chern class of the tautological quotient bundle \cite{3264}, and so
   \[ A^{2m-4}(\Gr(2,n))\ \xrightarrow{\cdot h}\ A^{2m-3}(\Gr(2,n)) \]
   is surjective (and hence an isomorphism, by hard Lefschetz). Let $\tau\colon Y\to \Gr(2,n)$ denote the inclusion morphism. The normal bundle formula tells us that the composition
   \[ A^{2m-4}(\Gr(2,n))    \ \xrightarrow{\cdot h}\ A^{2m-3}(\Gr(2,n)) \ \xrightarrow{\tau^\ast} \ A^{2m-3}(Y)\ \xrightarrow{\tau_\ast}\ A^{2m}(\Gr(2,n)) \]
   is a non-zero multiple of 
   \[ A^{2m-4}(\Gr(2,n))    \ \xrightarrow{\cdot h^4}\ A^{2m}(\Gr(2,n)) \ .\]
   This last map is the same as
   \[ H^{4m-8}(\Gr(2,n),\QQ)\  \ \xrightarrow{\cdot h^4}\ H^{4m}(\Gr(2,n),\QQ) \ ,\]    
   which is an isomorphism thanks to hard Lefschetz for the $(4m-4)$-dimensional variety $\Gr(2,n)$. This proves the required injectivity of $GDA^{2m-3}(Y)$ into cohomology.
   
   As for (\rom2), one can either prove this directly, or can reduce to (\rom1) via the generically defined isomorphism
    \[ A^1_{hom}(C)\ \xrightarrow{\cong}\ A^{2m-3}_{hom}(Y) \]
    given by Theorem \ref{main0}(\rom2).
    \end{proof}

 \begin{proposition}\label{Fr2} The following families have the Franchetta property:
 
 \noindent
 (\rom1) the family $\Cc\times_{B_0}\Cc\to B_0$;
 
 \noindent
 (\rom2) the family $\YY\times_B \YY\to B$.
 \end{proposition}
 
 \begin{proof} 
 
 \noindent
 (\rom1) Let $\bar{\Cc}\subset \Pf\times\bar{B}$ denote the projective closure of $\Cc$, and let us consider the projection
   \[ \pi\colon \bar{\Cc}\times_{\bar{B}}\bar{\Cc}\ \to\ \Pf\times\Pf\ .\]
   This is a {\em stratified projective bundle\/} (in the sense of \cite{FLV}). As such, \cite[Proposition 5.2]{FLV} implies the equality
      \begin{equation}\label{stra2} GDA^\ast_{B_0}(C\times C)= \ima\Bigl( A^\ast(\Pf^\circ\times \Pf^\circ)\to A^\ast(C\times C)\Bigr) +  \Delta_\ast GDA^\ast_{B_0}(C)\ ,\end{equation}
      where $\Delta\colon C\to C\times C$ is the inclusion along the diagonal, and $\Pf^\circ\subset\Pf$ denotes the non-singular locus of the Pfaffian.
       As $\Pf^\circ$ has the Chow--K\"unneth property \cite[Example 2.7]{GrPf}, $A^\ast(\Pf^\circ\times \Pf^\circ)$ is generated by $A^\ast(\Pf^\circ)\otimes A^\ast(\Pf^\circ)$.  The equality \eqref{stra2}
       thus simplifies to
       \begin{equation}\label{simpl} GDA^\ast_{B_0}(C\times C)= \Bigl\langle (p_j)^\ast \ima\bigl( A^j(\Pf^\circ)\to A^\ast(C)\bigr)\ ,\, \Delta_C\Bigr\rangle\ .\end{equation}
       
  We now proceed to check that $GDA^j_{B_0}(C\times C)$ injects into cohomology:
  
  In case $j=1$, we know that $\Delta_C$ is linearly independent from the decomposable classes 
  \[  \Bigl\langle (p_j)^\ast \ima\bigl(A^j(\Pf^\circ)\to A^\ast(C)\bigr)\Bigr\rangle    \]
  in cohomology (indeed, we may assume that $C$ has genus $>0$, for otherwise the statement is vacuously true). The required injectivity then reduces to Proposition \ref{Fr}(\rom2).
  
  In case $j=2$, we know that $A^1(\Pf^\circ)$ is 1-dimensional, generated by a hyperplane class $H$ (cf. Lemma \ref{1} below). Since $C\subset\PP^2$ is a plane curve,
  clearly we have an equality
    \[ \Delta_C\cdot (p_i)^\ast(H) = \sum_{r=0}^2 {1\over \deg C}\, (p_1)^\ast(H^r)\cdot (p_2)^\ast(H^{2-r})\ \ \ \hbox{in}\ A^2(C\times C)\ ,\]
 and so
  \[ \begin{split} GDA^2_{B_0}(C\times C)&=  \Bigl\langle (p_j)^\ast\ima\bigl( A^j(\Pf^\circ)\to A^\ast(C)\bigr)\ ,\, \Delta_C\Bigr\rangle \cap A^2(C\times C)\\
   &=    \Bigl\langle (p_j)^\ast \ima\bigl( A^j(\Pf^\circ)\to A^\ast(C)\bigr)\Bigr\rangle\cap A^2(C\times C)\ .\\
   \end{split}\]
  The required injectivity then reduces to Proposition \ref{Fr}(\rom2).
  
  In the above, we have used the following lemma:
  
  \begin{lemma}\label{1} Let $\Pf^\circ\subset\Pf$ denote (as above) the non-singular locus of the Pfaffian. We have
    \[ A^1(\Pf^\circ)\cong \QQ[H]\ .\]
      \end{lemma}
      
    \begin{proof}(of the lemma.) We consider
      \[ \wt{\Pf}:= \Bigl\{ (\omega,K)\in  \Pf\times \Gr(2,n)   \, \Big\vert\,  K\subset\ker\omega  \Bigr\} \ \ \ \subset\ \Pf\times \Gr(2,n)   \ . \]
    The projection $\wt{\Pf}\to\Gr(2,n)$ is a projective bundle (and so $\wt{\Pf}$ is smooth), and the projection $\wt{\Pf}\to\Pf$ is an isomorphism over the non-singular locus (and so
    $\wt{\Pf}\to\Pf$ is a resolution of singularities).       
    
    Being a projective bundle over a Grassmannian, $\wt{\Pf}$ has Picard number 2:
    \[ A^1(\wt{\Pf})=\QQ^2\ .\]
    The complement of (the isomorphic pre-image of) $\Pf^\circ$ inside $\wt{\Pf}$ is an irreducible divisor $D$ (it is a partial flag variety). The localization sequence 
    \[  A_\ast(D)\ \to\ A^1(\wt{\Pf})\ \to\ A^1(\Pf^\circ)\ \to\ 0 \]
    then gives the result.
     \end{proof}  
  
  \noindent
  (\rom2) Again, we may assume that $n$ is even (for otherwise the statement is vacuously fulfilled).
  In view of Lemma \ref{spread}, it will suffice to prove the Franchetta property for $\YY\times_{B_0} \YY\to B_0$.   Thanks to Theorem \ref{main0}(\rom2), for any fiber $Y=Y_b$ with $b\in B_0$ we have split injections
       \[ A^j(Y\times Y)\ \hookrightarrow\ A^{j+1-d}(C\times C)\oplus \bigoplus A^\ast(C)\oplus \QQ^s\ .\]
       The isomorphism of Theorem \ref{main0} being generically defined, there are also split injections
        \[ GDA^j(Y\times Y)\ \hookrightarrow\ GDA^{j+1-d}(C\times C)\oplus \bigoplus GDA^\ast(C)\oplus \QQ^s\ .\]  
        The required injectivity now follows from (\rom1) and Proposition \ref{Fr}(\rom2).
   \end{proof}

\subsection{MCK}

\begin{theorem}\label{main} Let $Y$ be a smooth dimensionally transverse intersection
   \[ Y:=\Gr(2,n)\cap H_1\cap H_2 \cap H_3 \ ,\]
   where the $H_j$ are Pl\"ucker hyperplanes.
 Then $Y$ has an MCK decomposition.
\end{theorem}

\begin{proof}
In case $n$ is odd, $Y$ has trivial Chow groups (Theorem \ref{main0}(\rom1)) and so the statement is vacuously true. In case $n=4$, $Y$ is a rational curve and again the statement is vacuously true. We may thus suppose that $n$ is even and $\ge 6$. We have the following general result:

\begin{proposition}\label{gen} Let $\YY\to B$ be a family of smooth projective varieties, verifying 

\noindent
(a1) the fibers $Y_b$ are of odd dimension $d\ge 5$ and 
  \[ A^j_{hom}(Y_b)=0\ \ \ \forall j> (d+1)/2\ \    \forall\ b\in B\ ;\]

\noindent
(a2) the fibers $Y_b$ have a generically defined K\"unneth decomposition, i.e. there exist $\{\pi^j_\YY\}\in A^d(\YY\times_B \YY)$ such that the fiberwise restriction
$\pi^j_{Y_b}:=\pi^j_\YY\vert_b \in A^d(Y_b\times Y_b)$ is a K\"unneth decomposition for all $b\in B$;

\noindent
(a3) the family $\YY\times_B \YY\to B$ has the Franchetta property.

Then $\{\pi^j_{Y_b}\}$ is an MCK decomposition for any $b\in B$.
\end{proposition}

\begin{proof}(of Proposition \ref{gen}.) Condition (a1) implies (via the Bloch--Srinivas argument, cf. \cite{BS}) that for every fiber $Y_b$ there exists a curve $C_b$ and a split injection
of motives
  \begin{equation}\label{inj} h(Y_b)\ \hookrightarrow\ h(C_b)((1-d)/2)\oplus \bigoplus \one(\ast)\ \ \ \hbox{in}\ \MM_{\rm rat}\ .\end{equation}
  
Condition (a3) implies that the K\"unneth decomposition $\{\pi^j_{Y_b}\}$ of (a2) is a self-dual CK decomposition. Let $h(Y_b)=\oplus_j h^j(Y_b)$ denote the corresponding decomposition of the motive of $X$. Using the injection \eqref{inj}, one finds that $h^j(Y_b)=\oplus \one(\ast)$ for all $j\not= d$, while for $j=d$ one finds a split injection
  \begin{equation}\label{inj2} h^d(Y_b)\ \hookrightarrow\ h^1(C_b)((1-d)/2)\ \ \ \hbox{in}\ \MM_{\rm rat}\ .\end{equation}
  
 Let us now establish that the CK decomposition $\{\pi^j_{Y_b}\}$ is MCK. By definition, what we need to check is that the cycle
   \[  \Gamma_{ijk}:=   \pi_{Y_b}^k\circ \Delta_{Y_b}^{\rm sm}\circ (\pi_{Y_b}^i\times \pi_{Y_b}^j)\ \ \ \in\ A^{2d}(Y_b\times Y_b\times Y_b) \]
   is zero for all $i+j\not=k$.
   
  Let us assume at least one of the integers $i,j,k$ is different from $d$. In this case, there is an injection
    \[ \Gamma_{ijk}\ \ \in \ ( \pi^{2d-i}_{Y_b}\times \pi^{2d-j}_{Y_b} \times \pi^k_{Y_b})_\ast  A^{2d}(Y_b\times Y_b\times Y_b)\ \hookrightarrow\ \bigoplus A^\ast(Y_b\times Y_b)\ ,\]
    and this injection sends generically defined cycles to generically defined cycles. But $\Gamma_{ijk}$ is generically defined and homologically trivial, and so the Franchetta property for $\YY\times_B \YY\to B$ gives the required vanishing  $\Gamma_{ijk}=0$.
    
   Next, let us assume $i=j=k=d$. In this case, the injection of motives \eqref{inj2}
    induces an injection of Chow groups
   \[ \Gamma_{ijk}\ \ \in\ ( \pi^{d}_{Y_b}\times \pi^{d}_{Y_b} \times \pi^{d}_{Y_b})_\ast  A^{2d}(Y_b\times Y_b\times Y_b)\ \hookrightarrow\   A^{ (d+3)/2 }(C_b\times C_b\times C_b)\ .\]
   But the right-hand side vanishes for dimension reasons for any $d\ge 5$, and so $\Gamma_{ijk}=0$.    
   \end{proof}

Let us now consider the family $\YY\to B$ of all smooth complete intersections $\Gr(2,n)\cap H_1\cap H_2 \cap H_3$, where $n\ge 6$ is even. Each fiber $Y_b$ has a generically defined CK decomposition $\{\pi^j_{Y_b}\}$ (Lemma \ref{ck}).
To check that $\{\pi^j_{Y_b}\}$ is MCK, it suffices to do this over a dense open of $B$; for instance we may take $B_0\subset B$ the locus as before where $Y_b$ has a smooth dual curve $C_b$ contained in $\Pf^\circ$.
Let us check that $\YY\to B_0$ verifies the conditions of Proposition \ref{gen}. Condition (a1) is immediate from Theorem \ref{main0}(\rom2). Condition (a2) is fulfilled by the $\{\pi^j_{Y_b}\}$.
As for condition (a3), this is Proposition \ref{Fr2}(\rom2). This ends the proof.
%
%
%
    \end{proof}

 \section{The tautological ring}
 
 \subsection{A positive result}
 
 \begin{corollary}\label{cor1} Let $Y$ be as in Theorem \ref{main}, and $m\in\NN$. Let
  \[ R^\ast(Y^m):=\Bigl\langle  (p_i)^\ast \ima\bigl( A^\ast(\Gr(2,n))\to A^\ast(Y)\bigr) , \, (p_{ij})^\ast(\Delta_Y)\Bigr\rangle\ \subset\ \ \ A^\ast(Y^m)   \]
  be the $\QQ$-subalgebra generated by (pullbacks of) cycles coming from $\Gr(2,n)$ and (pullbacks of) the diagonal $\Delta_Y\in A^d(Y\times Y)$. (Here $p_i$ and $p_{ij}$ denote the various projections from $Y^m$ to $Y$ resp. to $Y\times Y$).
  The cycle class map induces injections
   \[ R^\ast(Y^m)\ \hookrightarrow\ H^\ast(Y^m,\QQ)\ \ \ \hbox{for\ all\ }m\in\NN\ .\]
   \end{corollary}

\begin{proof} This is inspired by the analogous result for cubic hypersurfaces \cite[Section 2.3]{FLV3}, which in turn is inspired by analogous results for hyperelliptic curves \cite{Ta2}, \cite{Ta} (cf. Remark \ref{tava} below) and for K3 surfaces \cite{Yin}.

The Chow ring $A^\ast(\Gr(2,n))$ is generated by the Pl\"ucker polarization $h\in A^1(\Gr(2,n))$ and the Chern class $c_2(Q)\in A^2(\Gr(2,n))$, where $Q\to\Gr(2,n)$ is the universal quotient bundle \cite{3264}. As in \cite[Section 2.3]{FLV3}, let us write 
  \[  o:={1\over \deg (h^d)} \, h^d\ \ \in\  A^d(Y)\ , \ \ \ \ c:=c_2(Q)\vert_Y\ \ \in\ A^2(Y)\ ,\] 
 and
  \[ \tau:= \pi^d_Y=\Delta_Y - \, \sum_{j\not=d}  \pi^j_Y\ \ \in\ A^d(Y\times Y) \ ,\]
  where the $\pi^j_Y$ are as above, and $d:=\dim Y$.
  
Moreover, for any $1\le i<j\le m$ let us write 
  \[ \begin{split}   o_i&:= (p_i)^\ast(o)\ \ \in\ A^d(Y^m)\ ,\\
                        h_i&:=(p_i)^\ast(h)\ \ \in \ A^1(Y^m)\ ,\\
                        c_i&:=   (p_i)^\ast(c)\ \ \in\ A^2(Y^m)\ ,\\                       
                          \tau_{ij}&:=(p_{ij})^\ast(\tau)\ \ \in\ A^d(Y^m)\ .\\
                         \end{split}\]
 Note that (by definition) we have
   \[ R^\ast(Y^m)= \Bigl\langle    o_i, h_i, c_i, \tau_{ij}\Bigr\rangle\ \ \ \subset\ A^\ast(Y^m)\  .\]                   
                         
  Let us now define the $\QQ$-subalgebra
  \[ \bar{R}^\ast(Y^m):=\Bigl\langle o_i, h_i, c_i  , \tau_{ij}\Bigr\rangle \ \ \ \subset\   H^\ast(Y^m,\QQ) \]
  (where $i$ ranges over $1\le i\le m$, and $1\le i<j\le m$); this is the image of $R^\ast(Y^m)$ in cohomology. One can prove (just as \cite[Lemma 2.11]{FLV3} and \cite[Lemma 2.3]{Yin}) that the $\QQ$-algebra $ \bar{R}^\ast(Y^m)$
  is isomorphic to the free graded $\QQ$-algebra generated by $o_i,h_i,c_i, \tau_{ij}$, modulo the following relations:
    \begin{equation}\label{E:X'}
			\quad h_i \cdot o_i =c_i\cdot o_i= 0,   \quad c_i^{ (d+1)/2}=0,  \quad c_i^{(d-1)/2}=\lambda h_i^{d-1}   , \quad \ldots ,   \quad 
			h_i^d =\deg (h^d)\,o_i      \,;
			\end{equation}
			\begin{equation}\label{E:X2'}
			\tau_{ij} \cdot o_i =  \tau_{ij} \cdot h_i =\tau_{ij}\cdot c_i= 0, \quad \tau_{ij} \cdot \tau_{ij} = -b_d \, o_i\cdot o_j
			\,;
			\end{equation}
			\begin{equation}\label{E:X3'}
			\tau_{ij} \cdot \tau_{ik} = \tau_{jk} \cdot o_i\,;
			\end{equation}
			\begin{equation}\label{E:X4'}
			\sum_{\sigma \in \mathfrak{S}_{b_d+2}} 
			\prod_{i=1}^{b_d/2+1} \tau_{\sigma(2i-1), \sigma(2i)} = 0\, . 
			\end{equation}
Here $\lambda\in\QQ$, and the dots ``$\ldots$'' in \eqref{E:X'} indicate certain relations of type $c_i^{m_j} h_i^{n_j}=\lambda_j h_i^{2m_j +n_j}$.
By definition, $b_d:=\dim H^d(Y,\QQ)$ and $ \mathfrak{S}_{r}$ denotes the symmetric group on $r$ elements.

To prove Corollary \ref{cor1}, it suffices to check that all these relations  are verified modulo rational equivalence. 
The relations \eqref{E:X'} take place in $R^\ast(Y)$ and so they follow from the Franchetta property for $Y$ (Proposition \ref{Fr}). 
The relations \eqref{E:X2'} take place in $R^\ast(Y^2)$. The last relation is trivially verified, because ($Y$ being Fano) $A^{2d}(Y^2)=\QQ$. As for the other relations of \eqref{E:X2'}, these follow from the Franchetta property for $Y\times Y$ (Proposition \ref{Fr2}). 
   
 Relation \eqref{E:X3'} takes place in $R^\ast(Y^3)$ and follows from the MCK decomposition. Indeed, we have
   \[  \Delta_Y^{\rm sm}\circ (\pi^d_Y\times\pi^d_Y)=   \pi^{2d}_Y\circ \Delta_Y^{\rm sm}\circ (\pi^d_Y\times\pi^d_Y)  \ \ \ \hbox{in}\ A^{2d}(Y^3)\ ,\]
   which (using Lieberman's lemma) translates into
   \[ (\pi^d_Y\times \pi^d_Y\times\Delta_Y)_\ast    \Delta_Y^{\rm sm}  =   ( \pi^d_Y\times \pi^d_Y\times\pi^{2d}_Y)_\ast \Delta_Y^{\rm sm}                            
                                \ \ \ \hbox{in}\ A^{2d}(Y^3)\ ,\]
   which means that
   \[  \tau_{13}\cdot \tau_{23}= \tau_{12}\cdot o_3\ \ \ \hbox{in}\ A^{2d}(Y^3)\ .\]
   
  Finally, relation \eqref{E:X4'}, which takes place in $R^\ast(Y^{ b_d +2})$, is related to the Kimura finite-dimensionality relation \cite{Kim}:
  relation \eqref{E:X4'} expresses the vanishing
    \[ \sym^{b_d+2} H^{d}(Y,\QQ)=0\ ,\]
    where $H^{d}(Y,\QQ)$ is seen as a super vector space.
 This relation is also verified modulo rational equivalence, 
 (i.e., relation \eqref{E:X4'} is also true in $A^{d(b_d+2)}(Y^{b_d+2})$): relation \eqref{E:X4'} involves a cycle in
   \[ A^\ast(\sym^{b_d+2} h^d(Y))\ ,\]
  and $\sym^{b_d+2} h^d(Y)$ is $0$ 
 because $Y$ has Kimura finite-dimensional motive (Theorem \ref{main0}).

   This ends the proof.
 \end{proof}

\begin{remark}\label{tava} Given any curve $C$ and an integer $m\in\NN$, one can define the {\em tautological ring\/}
  \[ R^\ast(C^m):=  \bigl\langle  (p_i)^\ast(K_C),(p_{ij})^\ast(\Delta_C)\bigr\rangle\ \ \ \subset\ A^\ast(C^m) \]
  (where $p_i, p_{ij}$ denote the various projections from $C^m$ to $C$ resp. $C\times C$).
  Tavakol has proven \cite[Corollary 6.4]{Ta} that if $C$ is a hyperelliptic curve, the cycle class map induces injections
    \[  R^\ast(C^m)\ \hookrightarrow\ H^\ast(C^m,\QQ)\ \ \ \hbox{for\ all\ }m\in\NN\ .\]
    
   On the other hand, there exist curves for which the tautological ring $R^\ast(C^3)$ does {\em not\/} inject into cohomology, cf. Proposition \ref{neg} below.
\end{remark}

\subsection{A negative result}

 \begin{proposition}\label{neg} Let 
   \[ Y:=\Gr(2,n)\cap H_1\cap H_2\cap H_3 \]
   be a very general intersection of the Grassmannian with 3 Pl\"ucker hyperplanes, where $n$ is even and $8\le n\le  2000$. Let $C\subset \Pf$ be the curve dual to $Y$ (Definition \ref{dual}).
   The $\QQ$-subalgebra
     \[  R^\ast(C^m):=\bigl\langle (p_i)^\ast(K_C) ,\,  (p_{ij})^\ast(\Delta_C)\bigr\rangle\ \ \subset\ A^\ast(C^m) \]
    does {\em not\/} inject into cohomology for $m=3$.
  \end{proposition}
  
  \begin{proof} The point is that $C$ is a plane curve of degree $n/2$, and that the general plane curve of degree $n/2$ arises in this way \cite{Beau1}. 
  Using the spread lemma (Lemma \ref{spread}), it follows that the assumption that $R^\ast(C^m)$ injects into cohomology for the very general $C$ as in Proposition \ref{neg} would imply that $R^\ast(C^m)$ injects into cohomology
  for every plane curve of degree $n/2$. Taking $m=3$, this would mean that every plane curve of degree $n/2$ has a self-dual MCK decomposition. As explained in \cite[Proposition 7.1]{FV} and \cite[Remark 2.4]{FLV2},
  this would imply that for every plane curve $C$ of degree $n/2$ the Ceresa cycle 
    \[ C-[-1]_\ast(C)\ \ \in\  A_1(\hbox{Jac}(C))\] 
 is algebraically trivial.
   But this is known to be false for the Fermat curve of degree between 4 and 1000, cf. \cite{Ots}.
   \end{proof}

 \vskip0.5cm
\begin{nonumberingt} Thanks to Lie Fu and Charles Vial for lots of inspiring exchanges around MCK. 
\end{nonumberingt}

\end{document}